\documentclass[11pt,reqno]{amsart}
\usepackage{amssymb,amscd,amsbsy}
\usepackage{amsthm}
\newcommand{\lb}{\linebreak}

\let\pl = \l

\renewcommand{\a}{\alpha}
\renewcommand{\b}{\beta}
\newcommand{\g}{\gamma}
\renewcommand{\d}{\delta}

\renewcommand{\l}{\lambda}

\newcommand{\f}{\varphi}

\newcommand{\A}{{\mathcal A}}
\newcommand{\B}{{\mathcal B}}

\newcommand{\C}{{\Bbb C}}
\newcommand{\T}{{\Bbb T}}

\newcommand{\dd}{{\Bbb D}}
\newcommand{\R}{{\Bbb R}}
\newcommand{\Z}{{\Bbb Z}}

\newcommand{\rf}[1]{(\ref{#1})}

\newcommand{\const}{\operatorname{const}}

\newcommand{\eeq}{\end{equation}}
\newcommand{\beq}{\begin{equation}}
\newcommand{\bay}{\begin{eqnarray}}
\newcommand{\ba}{\begin{align*}}
\newcommand{\ea}{\end{align*}}
\newcommand{\ey}{\end{eqnarray}}
\newcommand{\bey}{\begin{eqnarray*}}
\newcommand{\eey}{\end{eqnarray*}}

\newcommand{\be}{\infty}

\newcommand{\bl}{\blacksquare}

\newcommand{\Pf}{{\bf Proof. }}

\newtheorem{thm}{\hspace{\parindent}Theorem}[section]

\newtheorem{lem}[thm]{\hspace{\parindent}Lemma}

\theoremstyle{remark}

\newtheorem*{rem*}{Remark}

\begin{document}

\newcommand{\vse}{\vspace{.2in}}
\numberwithin{equation}{section}

\title{\bf On S. Mazur's problems 8 and 88 \\from the Scottish Book}
\author{V.V. Peller}
\thanks{The author is partially supported by NSF grant DMS 0200712}
\maketitle

\begin{abstract}
The paper discusses Problems 8 and 88 posed by Stanis\pl aw Mazur in 
the Scottish Book \cite{SB}. It turns out that negative solutions to both
problems are immediate consequences of the results of \S 5 of \cite{P1}.
We discuss here some quantitative aspects of Problems 8 and 88 and give answers to open problems
discussed in a recent paper [PS] in connection with Problem 88.
\end{abstract}

\section{\bf Introduction}
\setcounter{equation}{0}

\

We are going to discuss in this paper Problems 8 and 88 posed by Stanis\pl aw Mazur in the Scottish Book
\cite{SB}. Problem 88 asks whether a Hankel matrix in the injective tensor product 
$\ell^1\check{\otimes}\ell^1$ of two spaces $\ell^1$
must have finite sum of the moduli of its matrix entries. Problem 8 asks whether for an arbitrary sequence
$\{z_n\}_{n\ge0}$ in the space $c$ of converging sequences
there exist sequences $\{x_n\}_{n\ge0}$ and $\{y_n\}_{n\ge0}$
in the space $c$ such that 
$$
z_n=\frac1{n+1}\sum_{k=0}^nx_ky_{n-k},\quad n\ge0.
$$
We give precise statements of the problems and all necessary definitions later.

It turned out that both problems have negative solutions. Independently, solutions were obtained by
Kwapie\'n and Pe\pl czy\'nski \cite{KP} and Eggermont and Leung \cite{EL}. In a recent paper by
Pe\pl czy\'nski and Sukochev \cite{PS} in Section 6 certain quantitative results related to negative 
solutions of Problems 8 and 88 are obtained and certain open problems are raised.

It turns out, however, that the results of Section 5 of my earlier paper \cite{P1} immediately
imply negative solutions to Problems 8 and 88. Moreover, Section 5 of \cite{P1} contains much stronger results.
In particular, a complete description of the Hankel matrices\footnote{Note that Hankel
matrices and Hankel operators play an important role in many areas of mathematics and applications, see
\cite{P2}.} in the injective tensor product of two 
spaces $\ell ^1$ is obtained in \cite{P1} in terms of the Besov space $B^1_{\be,1}$. 
Unfortunately, I was not aware about the Problems 8 and 88 when I wrote the paper \cite{P1}. 

In Sections 3 and 4 of this paper we explain why the results of \cite{P1} immediately imply negative
solutions to Problems 8 and 88 and we give a solution to the problems raised in \cite{PS}. 

In \S 2 we collect  necessary information on tensor products and Besov spaces.

\

\section{\bf Preliminaries}
\setcounter{equation}{0}

\

{\bf 1. Projective and injective tensor products.} We define the {\it projective tensor product}
$\ell^\be\hat{\otimes}\ell^\be$ as the space of matrices $\{q_{jk}\}_{j,k\ge0}$ of the form
\bay
\label{lbe}
q_{jk}=\sum_{n\ge0}a^{(n)}_jb^{(n)}_k,
\ey
where the $a^{(n)}=\{a^{(n)}_j\}_{j\ge0}$ and $b^{(n)}=\{b^{(n)}_j\}_{j\ge0}$ are sequences in
$\ell^\be$ such that
\bay
\label{fin}
\sum_{n\ge0}\|a^{(n)}\|_{\ell^\be}\|b^{(n)}\|_{\ell^\be}<\be.
\ey
The norm of the matrix $\{q_{jk}\}_{j,k\ge0}$ in $\ell^\be\hat{\otimes}\ell^\be$ is defined 
as the infimum of the left-hand side of \rf{fin} over all sequences $a^{(n)}=\{a^{(n)}_j\}_{j\ge0}$ and
$b^{(n)}=\{b^{(n)}_j\}_{j\ge0}$ satisfying \rf{lbe}.

Similarly, one can define the projective tensor products $c\hat{\otimes}c$ and $c_0\hat{\otimes}c_0$, where
$c$ is the subspace of $\ell^\be$ that consists of the converging sequences and $c_0$ is the subspace of $c$
that consists of the sequences with zero limit.

We consider the space $V^2$ that is a kind of a ``weak completion'' of $\ell^\be\hat{\otimes}\ell^\be$. $V^2$
consists of the matrices $Q=\{q_{jk}\}_{j,k\ge0}$ for which 
$$
\sup_{n>0}\|P_nQ\|_{\ell^\be\hat{\otimes}\ell^\be}<\be,
$$
where the projections $P_n$ are defined by
$$
(P_nQ)_{jk}=\left\{\begin{array}{ll}q_{jk},&j\le n,~k\le n,\\[.2cm]0,&\mbox{otherwise}.\end{array}\right.
$$
Note that
$c\hat\otimes c\subset \ell^\infty\hat\otimes\ell^\infty\subset V^2$.

The {\it injective tensor product} $\ell^1\check{\otimes}\ell^1$ of two spaces $\ell^1$ is, by definition, the
space of matrices $Q=\{q_{jk}\}_{j,k\ge0}$ such that
$$
\|Q\|_{\ell^1\check{\otimes}\ell^1}=\sup\left|\sum^N_{j,k=0}q_{jk}x_jy_k\right|<\be,
$$
where the supremum is taken over all sequences $\{x_j\}_{j\ge0}$ and $\{y_k\}_{k\ge0}$ in the unit ball of
$\ell^\be$ and over all positive integers $N$. The space $\ell^1\check{\otimes}\ell^1$ can naturally be
identified with the space of bounded linear operators from $c_0$ to $\ell^1$ (note that every bounded operator
from $c_0$ to $\ell^1$ is compact).

\medskip

{\bf2. Besov spaces.} In this paper we consider only Besov spaces of functions analytic in the unit disk $\dd$.
Besov spaces $B_{p,q}^s$ admit many equivalent descriptions. We give a definition
in terms of dyadic Fourier expansions. We define the polynomials $W_n$, $n\ge0$, as follows. If $n\ge1$, then
$\widehat W_n(2^n)=1$, $\widehat W_n(k)=0$ for $k\not\in(2^{n-1},2^{n+1})$, and $\widehat W_n$ is a linear
function on $[2^{n-1},2^n]$ and on $[2^n,2^{n+1}]$. We put $W_0(z)=1+z$. It is easy to see that
$$
\|W_n\|_{L^1}\le\frac32,\quad n\ge0,
$$
and
$$
f=\sum_{n\ge0}f*W_n
$$
for an arbitrary analytic function $f$ in $\dd$.

For $1\le p\le\be$, $1\le q\le\be$, and $s\in\R$, we define the Besov space $B^s_{p,q}$ as the space of
analytic functions in $\dd$ satisfying
\bay
\label{bes}
f\in B^s_{p,q}\quad\Longleftrightarrow\quad
\{2^{ns}\|f*W_n\|_{L^p}\}_{n\ge0}\in\ell^q.
\ey
If $q=\be$, the space $B^s_{p,q}$ is nonseparable. We denote by $b^s_{p,\be}$ the closure of the set of
polynomials in $B^s_{p,\be}$. It is easy to verify that
$$
f\in b^s_{p,\be}\quad\Longleftrightarrow\quad
\{2^{ns}\|f*W_n\|_{L^p}\}_{n\ge0}\in c_0.
$$

Besov spaces admit many other descriptions (see \cite{Pe} and \cite{P2}).

\

\section{\bf Problem 8}
\setcounter{equation}{0}

\

To state Mazur's Problem 8 of the Scottish Book \cite{SB}, consider the bilinear form
$\B$ on $c\times c$ defined by
$$
\B\big(\{x_n\}_{n\ge0},\{y_n\}_{n\ge0}\big)=\{z_n\}_{n\ge0},
$$
where
$$
z_n=\frac1{n+1}\sum_{k=0}^nx_ky_{n-k},\quad n\ge0,
$$
and $c$ is the space of sequences that have a limit at infinity.

It is easy to see that $B$ maps $c\times c$ into $c$. {\it S. Mazur asked in Problem 8 whether 
$\B$ maps $c\times c$ onto $c$}.

As mentioned in the Introduction, a negative solution to problem 8 follows immediately from
the results of \S 5 of \cite{P1}. To state Theorem 5.1 of \cite{P1}, we define the operator $\A$ on the space of
matrices. Let $Q=\{q_{jk}\}_{j,k\ge0}$. Then $\A Q$ is the sequence defined by
$$
\A Q=\{z_n\}_{n\ge0},\quad\mbox{where}\quad z_n=\frac1{n+1}\sum_{j+k=n}q_{jk}.
$$

\medskip

{\bf Theorem 5.1 of [P1].} {\it$\A$ maps the space $V^2$ onto the space of Fourier coefficients of the
Besov space $B^0_{1,\be}$.}

\medskip

Recall that the space $V^2$ is defined in the introduction. In particular, it follows from Theorem 5.1 of [P1]
that 
$$
\A(c\hat{\otimes}c)\subset\A(\ell^\be\hat{\otimes}\ell^\be)\subset\A(V^2)\subset
\left\{\{\hat f(n)\}_{n\ge0}:~f\in B^0_{1,\be}\right\},
$$
and so
$$
\B(c\times c)\subset\left\{\{\hat f(n)\}_{n\ge0}:~f\in B^0_{1,\be}\right\}.
$$
It is easy to see that
$$
c\not\subset\left\{\{\hat f(n)\}_{n\ge0}:~f\in B^0_{1,\be}\right\}.
$$
Indeed, if $f\in B^0_{1,\be}$, then it follows immediately from \rf{bes} and from \cite{R}, \S 8.6 that
$$
\sup_{n\ge0}\sum_{k=0}^n\big|\hat f(2^n+2^k)\big|^2<\be.
$$
This gives a negative solution to Problem 8.

In fact, Theorem 5.1 of \cite{P1} allows one to describe $\A(c\hat{\otimes}c)$. First, let us
observe that Theorem 5.1 of \cite{P1} easily implies the following description of $\A(c_0\hat{\otimes}c_0)$.

\begin{thm}
\label{c0}
$$
\A(c_0\hat{\otimes}c_0)=\left\{\{\hat f(n)\}_{n\ge0}:~f\in b^0_{1,\be}\right\}.
$$
\end{thm}

Recall that $b^0_{1,\be}$ is the closure of the polynomials in $B^0_{1,\be}$ (see \S 2). Theorem \ref{c0},
in turn, easily implies the following description of $\A(c\hat{\otimes}c)$.

\begin{thm}
\label{c}
$$
\A(c\hat{\otimes}c)=\left\{\{\hat f(n)+d\}_{n\ge0}:~f\in b^0_{1,\be},~d\in\C\right\}.
$$
\end{thm}

\

\section{\bf Problem 88}
\setcounter{equation}{0}

\

Recall that in Problem 88 of \cite{SB} {\it S. Mazur asked whether a Hankel matrix
$\{\g_{j+k}\}_{j,k\ge0}$ in the injective tensor product $\ell^1\check{\otimes}\ell^1$
must satisfy the condition}:
$$
\sum_{k=0}^\be(1+k)|\g_k|<\be,
$$
i.e., whether the sum of the moduli of the matrix entries must be finite.

As mentioned in the Introduction, a negative solution to Problem 88 follows immediately from the results
of \S 5 of \cite{P1}.  A complete description of Hankel matrices in $\ell^1\check{\otimes}\ell^1$
is given by Theorem 5.2 of \cite{P1}:

\medskip

{\bf Theorem 5.2 of [P1].} {\it A Hankel matrix $\{\g_{j+k}\}_{j,k\ge0}$ belongs to 
$\ell^1\check{\otimes}\ell^1$ if and only if the function $f$ defined by
$$
f(z)=\sum_{n\ge0}\g_n z^n
$$
belongs to the Besov class $B^1_{\be,1}$.}

\medskip

Let us obtain the best possible estimate on the moduli of the matrix entries of Hankel matrices in 
$\ell^1\check{\otimes}\ell^1$.

Since $\|f*W_n\|_{L^2}\le\|f\|_{L^2}\|W_n\|_{L^1}\le 3/2\|f\|_{L^2}$,
it follows easily from \rf{bes} that if $f\in B^1_{\be1}$, then 
\bay
\label{mod}
\sum_{n=0}^\infty2^n\left(\sum_{k=2^n}^{2^{n+1}-1}|\hat f(k)|^2\right)^{1/2}<\infty.
\ey
Let us show that this is the best possible estimate for the moduli of the Fourier coefficients
of functions in $B^1_{\be1}$. To show this, we are going to use a version of the 
de Leeuw--Katznelson--Kahane theorem. It was proved in \cite{dLKK} that if $\{\b_n\}_{n\in\Z}$
is a sequence of nonnegative numbers in $\ell^2(\Z)$, then there exists a continuous function $f$ on $\T$ such
that
$$
|\hat f(n)|\ge\b_n,\quad n\in\Z.
$$
We refer the reader to \cite{K1}, \cite{K2}, and \cite{N} for refinements of the 
de Leeuw--Katznelson--Kahane theorem and different proofs. We need the following version of the
de Leeuw--Katznelson--Kahane theorem: 

\begin{lem}
\label{kis}
There is a positive number $K$ such that for arbitrary nonnegative numbers
$\b_0,\b_1,\cdots,\b_m$, there
exists a polynomial $f$ of degree $m$ such that
$$
|\hat f(j)|\ge\b_j,\quad 0\le j\le n,\quad\mbox{and}\quad\|f\|_{L^\be}\le
K\left(\sum_{j=0}^n\b_j^2\right)^{1/2}.
$$
\end{lem}

Lemma \ref{kis} follows easily from the results of \cite{K2}.

\begin{thm}
\label{LKK}
Let $\{\a_k\}_{k\ge0}$ be a sequence of nonnegaitive numbers such that 
\bay
\label{maj}
\sum_{n=0}^\infty2^n\left(\sum_{k=2^n}^{2^{n+1}-1}\a_k^2\right)^{1/2}<\infty.
\ey
Then there exists a function $\f$ in the space $B_{\be1}^1$ such that
 $|\hat\f(k)|\ge \a_k$ for $k\ge0$.
\end{thm}

\Pf By Lemma \ref{kis}, there exists $K>0$ and a sequence of polynomials $f_n$, $n\ge0$, such that
$$
f_0(z)=\hat f_0(0)+\hat f_0(1)z,\quad f_n(z)=\sum_{k=2^n}^{2^{n+1}-1}\hat f_n(k)z^k,
\quad\mbox{for}\quad n\ge1,
$$ 
$$
|\hat f_0(k)|\ge\a_k,\quad\mbox{for}\quad k=0,\,1,\quad|\hat f_n(k)|\ge\a_k\quad\mbox{for}
\quad n\ge1,~~2^n\le k\le2^{n+1}-1,
$$
and
$$
\|f_0\|_{L^\be}\le K(\a_0^2+\a_1^2)^{1/2},
\quad\|f_n\|_{L^\be}\le K\left(\sum_{k=2^n}^{2^{n+1}-1}\a_k^2\right)^{1/2}
\quad\mbox{for}\quad n\ge1.
$$

We can define now the function $\f$ by
$$
\f=\sum_{n\ge0}f_n.
$$
Obviously, $|\hat\f(k)|\ge \a_k$ for $k\ge0$. Let us show that $\f\in B_{\be1}^1$. We have
\begin{align*}
\qquad\sum_{n\ge1}2^n\|\f*W_n\|_{L^\be}&=\sum_{n\ge1}2^n\|\big(f_{n-1}+f_n+f_{n+1}\big)*W_n\|_{L^\be}\\[.2cm]
&\le\sum_{n\ge1}2^n\|\big(f_{n-1}+f_n+f_{n+1}\big)\|_{L^\be}\|W_n\|_{L^1}\\[.2cm]
&\le3\sum_{n\ge1}2^n\|f_n\|_{L^\be}\|W_n\|_{L^1}
\le3\cdot\frac32\sum_{n\ge1}2^n\|f_n\|_{L^\be}\\[.2cm]
&\le\frac92K
\sum_{n=0}^\infty2^n\left(\sum_{k=2^n}^{2^{n+1}-1}\a_k^2\right)^{1/2}<\be.\qquad\bl
\end{align*}

In \cite{PS} the following problem was considered. Let $\Psi$ be the function on $(0,\be)$ defined by
$$
\Psi(t)=\left\{\begin{array}{ll}\frac32t-1,&0<t\le2,\\[.2cm]
t,& t>2.\end{array}\right.
$$
Let $\{\g_{j+k}\}_{j,k\ge0}$ be a Hankel matrix. The following result was proved
in \cite{PS} (Theorem 6.7):

(i) {\em if $\b<\Psi(t)$, then}
$$
\sum_{k\ge0}|\g_k|^t(1+k)^\b<\be\quad\mbox{whenever}\quad\{\g_{j+k}\}_{j,k\ge0}\in\ell^1\check{\otimes}\ell^1;
$$

(ii) {\em if $\b>\Psi(t)$, then}
$$
\sum_{k\ge0}|\g_k|^t(1+k)^\b=\be\quad\mbox{for some}\quad\{\g_{j+k}\}_{j,k\ge0}\in\ell^1\check{\otimes}\ell^1;
$$

(iii) {\em if $\b=\Psi(t)$ and $\frac43\le t<\be$, then}
$$
\sum_{k\ge0}|\g_k|^t(1+k)^\b<\be\quad\mbox{whenever}\quad\{\g_{j+k}\}_{j,k\ge0}\in\ell^1\check{\otimes}\ell^1.
$$

In \cite{PS} the problem is raised to find out whether
$$
\sum_{k\ge0}|\g_k|^t(1+k)^{\Psi(t)}
$$
has to be finite for $t\in\big(0,\frac43\big)$ whenever
$\{\g_{j+k}\}_{j,k\ge0}\in\ell^1\check{\otimes}\ell^1$.

It is easy to deduce Theorem 6.7 of \cite{PS} from \rf{mod} and above Theorem \ref{LKK}. Moreover,
using \rf{mod} and Theorem \ref{LKK}, we can solve the problem posed in \cite{PS} and settle the case 
$t\in\big(0,\frac43\big)$.

\begin{thm}
\label{re}
If $1\le t<\frac43$, then 
$$
\sum_{k\ge0}|\g_k|^t(1+k)^{3t/2-1}<\be
\quad\mbox{whenever}\quad\{\g_{j+k}\}_{j,k\ge0}\in\ell^1\check{\otimes}\ell^1.
$$
If $0<t<1$, then 
$$
\sum_{k\ge0}|\g_k|^t(1+k)^{3t/2-1}=\be
\quad\mbox{for some}\quad\{\g_{j+k}\}_{j,k\ge0}\in\ell^1\check{\otimes}\ell^1.
$$
\end{thm}

\Pf Suppose that $1\le t<2$. By H\"older's inequality, we have
\begin{align*}
\sum_{k\ge1}|\g_k|^t(1+k)^{3t/2-1}&\le\,
\const\sum_{n\ge0}2^{n(3t/2-1)}\sum_{k=2^n}^{2^{n+1}}|\g_k|^t\\[.2cm]
&\le\sum_{n\ge0}2^{\frac{3n}{2}t}2^{-n}\left(\sum_{k=2^n}^{2^{n+1}}|\g_k|^2\right)^{t/2}
2^{n(1-t/2)}\\[.2cm]
&=\sum_{n\ge0}2^{nt}\left(\sum_{k=2^n}^{2^{n+1}}|\g_k|^2\right)^{t/2}.
\end{align*}
Since $t\ge1$, the $\ell^t$ norm of a sequence does not exceed its $\ell^1$ norm, and so
$$
\sum_{n\ge0}2^{nt}\left(\sum_{k=2^n}^{2^{n+1}}|\g_k|^2\right)^{t/2}\le
\left(\sum_{n\ge0}2^n\left(\sum_{k=2^n}^{2^{n+1}}|\g_k|^2\right)^{1/2}\right)^t.
$$
The result follows now from Theorem 5.2 of \cite{P1} and \rf{mod}.

Suppose now that $0<t<1$. It follows from Theorem \ref{LKK} that it suffices to find a sequence
$\{\a_k\}_{k\ge0}$ of nonnegative numbers that satisfies \rf{maj} and such that
$$
\sum_{k\ge0}\a_k^t(1+k)^{3t/2-1}=\be.
$$
Let $\{\d_n\}_{n\ge0}$ be a sequence of positive numbers such that $\{2^{3n/2}\d_n\}_{n\ge0}\in\ell^1$, but 
\lb$\{2^{3n/2}\d_n\}_{n\ge0}\not\in\ell^t$.

Put
$$
\a_0=0 \quad\mbox{and}\quad\a_k=\d_n\quad\mbox{if}\quad2^n\le k\le2^{n+1}-1.
$$
We have
$$
\sum_{n\ge0}2^n\left(\sum_{k=2^n}^{2^{n+1}-1}\a_k^2\right)^{1/2}=
\sum_{n\ge0}2^{3n/2}\d_n<\be.
$$
However,
\begin{align*}
\qquad\sum_{k\ge0}\a_k^t(1+k)^{3t/2-1}&\ge\const\sum_{n\ge0}2^{n(3t/2-1)}\sum_{k=2^n}^{2^{n+1}}\a_k^t\\[.2cm]
&=\const\sum_{n\ge0}2^{n(3t/2-1)}2^n\d_n^t=\const\sum_{n\ge0}2^{3nt/2}\d_n^t=\be.\quad\bl
\end{align*}

\

\

\noindent
\begin{tabular}{p{8cm}p{14cm}}
Department of Mathematics \\
Michigan State University  \\
East Lansing, Michigan 48824\\
USA\\
email: peller@math.msu.edu
\end{tabular}

\end{document}